\documentclass[a4paper,12pt]{article}
\usepackage{times, url}
\usepackage{amsmath,amssymb,amsthm,amsfonts,bm,cite,soul}

\newtheorem{theorem}{Theorem}[section]
\newtheorem{proposition}{Proposition}[section]
\newtheorem{corollary}{Corollary}[section]
\newtheorem{lemma}{Lemma}[section]
\newtheorem{remark}{Remark}[section]

\usepackage[center]{caption}
\usepackage[none]{hyphenat}

\usepackage{color, colortbl}
\definecolor{Gray}{gray}{0.9}

\begin{document}
	\setcounter{page}{1}

	\begin{center}
		{\LARGE \bf   Powers of the operator \boldmath$u(z)\frac{d}{dz}$ and their connection 
			\vspace{3pt}
			
			with some combinatorial numbers}
		\vspace{8mm}
		
		{\Large \bf Ioana Petkova}
		\vspace{3mm}
		
		Faculty of Mathematics and Informatics, Sofia University	\\
		Sofia, Bulgaria \\
		e-mail: \url{jepetkova@fmi.uni-sofia.bg}, 
		
		\ \url{ioana_petkovaa@abv.bg},

	\end{center}
	\vspace{10mm}
	
	\noindent
	{\bf Abstract:} 
	In this paper the operator $A = u(z)\frac{d}{dz}$ is considered, where $u$ is an entire or meromorphic function in the complex plane. The expansion of $A^{k}$ ($k\geq1$) with the help of the powers of the differential operator $D=\frac{d}{dz}$ is obtained, and it is shown that this expansion depends on special numbers. Connections between these numbers and known combinatorial numbers are given. Some special cases of the operator $A$, corresponding to $u(z) = z$, $u(z) = e^{z}$, $u(z) = \frac{1}{z}$, are considered. 
	
	
	\vspace{5mm}
	
	\section{Introduction} \label{sec:Intr}
	
	Let $ u: \mathbb{C} \rightarrow \mathbb{C}$ be an entire or meromorphic function. We denote by A the operator
	$$A = u(z)\frac{d}{dz}.$$
	We shall consider the operator $A$ over the set $U$, which is the union of the set of all entire functions in $\mathbb{C}$ and the set of all meromorphic functions in $\mathbb{C}$.
	Particular cases of $A$ are the operators: $A_{1} = z\frac{d}{dz}$,  $A_{2} = e^{z}\frac{d}{dz}$, $A_{3} = \frac{1}{z}\frac{d}{dz}$. 
	
	The operator $A_{1}$ is well studied. In \cite{Ross}, p.~242, $A_{1}$ is called Cauchy-Euler operator, although Euler first considered this operator. Thus we may call the operator $A$ generalized Euler-Cauchy operator. In \cite{Watson}, $A_{1}$ is called by Watson theta operator and is denoted by $\theta$. In \cite{7}, p.~542, $A_{1}$ is denoted by $\delta$ and it is shown that $\delta(\delta - 1) \dots (\delta -n+1) = z^{n} D^{n}$. 
	
	The operator $A_{3}$ is connected with the Bessel functions of the first kind $J=J_{\nu}(z)$, which satisfy the differential equation 
	
	\[
	z^{2} \left(\frac{d}{dz}\right)^{2}J + z \frac{d}{dz}J + (z^{2} - \nu^{2})J = 0. 
	\]
	In particular, for $\nu = n + \frac{1}{2}$ and $\nu = -n - \frac{1}{2}$ $(n = 0, 1, 2, \dots)$, we have: 
	\begin{align*}
		J_{n+\frac{1}{2}} (z) &= \sqrt{\frac{2}{\pi}}z^{n+\frac{1}{2}}(-1)^{n}A_{3}^{n} \frac{\text{sin}z}{z}; 	\\[2pt]
		J_{-n-\frac{1}{2}} (z) &= \sqrt{\frac{2}{\pi}}z^{n+\frac{1}{2}}A_{3}^{n} \frac{\text{cos}z}{z}
	\end{align*}	
	(see \cite{Kerimov} and \cite{JEL}, p. 178, p. 198).
	
	The operators $A_{1}\!^{k}$ and  $A_{3}\!^{k}$ belong to the class of the so-called hyper-Bessel differential operators  (Bessel type operators), i.e. operators of the form 
	
	\[
	B = z^{\alpha_{0}}\frac{d}{dz} z^{\alpha_{1}}\frac{d}{dz} \dots z^{\alpha_{m-1}}\frac{d}{dz}z^{\alpha_{m}},
	\]
	where $\alpha_{i}$, $i=0, \dots m$, are arbitrary real numbers such that $m > \displaystyle \sum_{i=0}^{m} \alpha_{i}$ \; (see \cite{Vi}, pp. 99-101).  
	
	Using the eigen functions of the operators $A_{1}$, $A_{2}$, $A_{3}$, one may obtain the following results:

	\begin{theorem}\label{thm:theorem1.1}
		Let $n$ be a natural number and  $P(w) = c_{0} + c_{1} w + c_{2} w^{2} + \dots + c_{n} w^{n}$, where $c_{t} \in \mathbb{C}$, $t = 0, 1, \dots, n$. Then for every $\lambda \in \mathbb{C}$ the equalities: 
	\end{theorem}
	
	\begin{equation*}
		P(A_{1}) z^{\lambda} = P(\lambda) z^{\lambda};
	\end{equation*}
	
	\begin{equation*}
		P(A_{2}) e^{-\lambda e^{-z}} = P(\lambda) e^{-\lambda e^{-z}};
	\end{equation*}
	
	\begin{equation*}
		P(A_{3}) e^{\frac{1}{2}\lambda z^{2}} = P(\lambda) e^{\frac{1}{2}\lambda z^{2}},
	\end{equation*}
	hold.
	
	\begin{proof}
		The first equality is proved in	\cite{Polya}, problem 44, p. 8. For the second equality, we verify that 
		
		\begin{align*}
			A_{2} e^{-\lambda e^{-z}} &= \lambda e^{-\lambda e^{-z}}, \\[-8pt]
			\intertext{ hence } \nonumber\\[-16pt]
			A_{2}\!^{t} \ e^{-\lambda e^{-z}} &= \lambda ^{t} e^{-\lambda e^{-z}}, \quad t=0, 1, 2, 3, ...
		\end{align*}	
		and the second equality holds. 
		
		For the third equality, we verify that 
		
		\begin{align*}
			A_{3}  e^{\frac{1}{2}\lambda z^{2}} &= \lambda e^{\frac{1}{2}\lambda z^{2}}, \\[-8pt]
			\intertext{hence} \\[-16pt]
			A_{3} \!^{t} \ e^{\frac{1}{2}\lambda z^{2}} &= \lambda ^{t} e^{\frac{1}{2}\lambda z^{2}}, \quad t=0, 1, 2, 3, ... 
		\end{align*}
		and the third equality holds. 
	\end{proof}
	
	As a corollary of Theorem \ref{thm:theorem1.1}
	, we obtain: 
	
	\begin{proposition}\label{prop:proposition1.1}
		For every two complex constants $a$ and $b$, the following equalities hold:
		
		\begin{equation*}
			P(a + bA_{1}) z^{\lambda} = P(a+b\lambda) z^{\lambda} \quad  \text{(see \cite{Belbahri},  p.551)}; 
		\end{equation*}
		
		\begin{equation*}
			P(a+bA_{2}) e^{-\lambda e^{-z}} = P(a+b\lambda) e^{-\lambda e^{-z}};
		\end{equation*}
		
		\begin{equation*}
			P(a+bA_{3}) e^{\frac{1}{2}\lambda z^{2}} = P(a+b\lambda) e^{\frac{1}{2}\lambda z^{2}}.
		\end{equation*}
	\end{proposition}
	
	\hfill $\square$
	
	In the same way, one may establish the following result:
	
	\begin{theorem}\label{thm:theorem1.2}
		Let $\lambda \in \mathbb{C}$ and $y=y(z)$ satisfy the equality $Ay = \lambda y$, where $A$ is the operator $A = u(z) \frac{d}{dz}$. Then for every $a, b \in \mathbb{C}$ the equality 
		
		\[
		P(a+bA) y = P(a+b\lambda)y
		\]
		holds.
	\end{theorem}
	
	\begin{remark}\label{rem:remark1.1}
		From the condition of the theorem, the differential equation $u(z) y' (z) = \lambda y (z)$ holds and its solution is given by
		
		\[
		y(z) = \tilde{c} \ e ^{\lambda \int_{c}^{z} \frac{dt}{u(t)}}
		,\]
		where $c, \tilde{c} \in \mathbb{C}$ are arbitrary constants.
		Let $\tilde{c} =1$. When $A = A_{1}$, we set $c=1$ and obtain $y(z) = z^{\lambda}$. When $A = A_{2}$, we set $c=\infty$ and obtain $y(z) = e^{-\lambda e^{-z}}$. When $A = A_{3}$, we set $c=0$ and obtain $y(z) = e^{\frac{1}{2}\lambda z^{2}}$.
	\end{remark}
	
	\bigskip
	
	Let $S(k,s)$ and $\sigma(k,s)$ be the Stirling numbers of the first and second kind, respectively (for the definition and properties of these numbers see \cite{Abr}, p.~824, p.~825). Then the following expansions are valid:  
	
	\begin{align}
		A_{1} \! ^{k} &=  \sum_{s=0}^{k}\sigma(k,s)z^{s}\left(\frac{d}{dz}\right)^{s}
		\text{	\quad   (see \cite{Riordan}, p. 218 and \cite{7}, p.~543)};\label{a}\\
		A_{2} \! ^{k} &= e^{kz}\sum_{s=0}^{k}(-1)^{k-s}S(k,s)\left(\frac{d}{dz}\right)^{s}; \label{b}\\
		A_{3} \! ^{k} &= \sum_{s=0}^{k}(-1)^{k-s}(2k-2s-1)!!{2k-s-1 \choose s-1}\frac{1}{z^{2k-s}}\left(\frac{d}{dz}\right)^{s},\label{c}
	\end{align}
	
	\noindent where we define (-1)!!=1.
	
	We shall prove the last two representations in the Appendix. Here we note that one may use the first two of the above equalities as defining equalities for introducing the Stirling numbers. 
	
	Just as every infinite sequence of numbers $\displaystyle\{a_{n}\}_{n=0}^{\infty}$ is connected with the generating function $f(z) = \displaystyle\sum_{n=0}^{\infty}a_{n}z^{n}$, it is interesting to know if the powers $A^{k}$ of $A$ are generating operators for some special sequences of numbers. 
	
	The mentioned above may motivate one to find an explicit representation for the operators $A^{k}$, $k=1,2,3,...$
	In this paper we show that the powers $A^{k}$ of $A$ are generating operators for different kinds of combinatorial numbers (in particular, for Stirling numbers of the first and second kind). 
	Also, several new identities for some combinatorial numbers are obtained.

	\section{The main result}
	
	First, we want to express the operator $A^{k}$ with the help of $\big(\frac{d}{dz}\big)^{s}$, $s=1, \dots, k$. 
	
	\begin{theorem}
		There exist numbers $C_{\alpha_{1}, ..., \alpha_{s}}^{s, m, k+1}$, for which the following representation:
		
		\begin{equation}\label{1}
			A^{k+1} = \sum_{s=1}^{k+1}P_{k+1}^{s}(u)\bigg(\frac{d}{dz}\bigg)^{s},
		\end{equation}
		is valid, where: 
		\begin{align}
			P_{k+1}^{k+1}(u) &= u^{k+1}; \label{2}\\
			P_{k+1}^{k+1-s}(u) &= \sum_{m=1}^{s}u^{k+1-m}F_{m}^{k+1}(u), s = 1, ..., k, \label{3}\\
			F_{m}^{k+1}(u) &= \sum_{m,s} C_{\alpha_{1}, ..., \alpha_{s}}^{s, m, k+1}\bigg(\frac{du}{dz}\bigg)^{\alpha_{1}}...\bigg(\frac{d^{s}u}{dz^{s}}\bigg)^{\alpha_{s}} \label{4}
		\end{align}
		and $\displaystyle\sum_{m,s}$ in (\ref{4}) is over all nonnegative integers $\alpha_{1}$, ..., $\alpha_{s}$ for which:
		\vspace{-15pt}
		\[
		\sum_{i=1}^{s}\alpha_{i} = m; \ \  \sum_{i=1}^{s} i\alpha_{i} = s. 
		\]
		The numbers  $C_{\alpha_{1}, ..., \alpha_{s}}^{s, m, k+1}$ satisfy the following recurrence relations:
		
		\small
		\begin{equation}	\label{5}
			C_{\gamma_{1}, ..., \gamma_{s}, \gamma_{s+1}}^{s+1, m, k+1} =	\sum_{t=s+1}^{k}\big[(t+1-m)C_{\gamma_{1} - 1, \gamma_{2},..., \gamma_{s}}^{s, m-1, t} +
			\sum_{j=1}^{s+1}(\gamma_{j}+1)C_{\gamma_{1}, ...,\gamma_{j-1}, \gamma_{j}+1, \gamma_{j+1}-1, \gamma_{j+2},..., \gamma_{s}}^{s, m, t}\big],
		\end{equation} 
		\normalsize
		where: $1 < m < s+1$; $\gamma_{s+1} = 0$; \ $\displaystyle\sum_{i=1}^{s+1}\gamma_{i} = m; \  \sum_{i=1}^{s+1}i\gamma_{i} = s+1; s = 1, ..., k-1; $

		\begin{equation}\label{6}
			C_{\gamma_{1}, ..., \gamma_{s}, \gamma_{s+1}}^{s+1, s+1, k+1} = \sum_{t=s+1}^{k}(t-s)C_{\gamma_{1} - 1, \gamma_{2}, ..., \gamma_{s}}^{s, s, t},
		\end{equation}
		where:
		$ \displaystyle
		\gamma_{s+1} = 0; \ \sum_{i=1}^{s+1}\gamma_{i} = s+1; \  \sum_{i=1}^{s+1}i\gamma_{i} = s+1; \ s=1, ..., k-1.
		$ 
		
		Also, the relation

		\begin{equation}\label{7}
			C_{\;\underbrace{\text{\scriptsize $0, ..., 0$}}_{\text{\normalfont s}}, 1}^{s+1, 1, k+1} = \sum_{t=s+1}^{k}C_{\;\underbrace{\text{\scriptsize $0, ..., 0$}}_{\text{\normalfont s-1}}, 1}^{s, 1, t}
		\end{equation}
		holds, where $s = 1, ..., k-1$. 
		$_{\square}$
	\end{theorem}
	
	\begin{remark}\label{rem:remark2.1}
		It is easy to see that (\ref{6}) takes the form: 
		\begin{equation}\label{8}
			C_{s+1, \;\underbrace{\text{\scriptsize $0, ..., 0$}}_{\text{\normalfont s}}}^{s+1, s+1, k+1} = \sum_{t=s+1}^{k}(t-s)C_{s,\;\underbrace{\text{\scriptsize $0, ..., 0$}}_{\text{\normalfont s-1}}}^{s, s, t}, \; s = 1, ..., k-1. 
		\end{equation}
	\end{remark}
	
	\begin{proof}[Proof of Theorem 3.] Let $f \in U$ and 
		\[
		A^{k}f(z) = \sum_{s=1}^{k}P_{k}^{s}(u)f^{(s)}(z).
		\]
		Since $A^{n}$, $n \geq1$, is a linear operator, 
		then 
		\begin{align*}
			&\sum_{s=1}^{k+1}P_{k+1}^{s}(u)f^{(s)}(z)
			= 	A^{k+1}f(z)  = A\big(A^{k}f(z)\big) = A \bigg(\sum_{s=1}^{k}P_{k}^{s}(u)f^{(s)}(z)\bigg) = \\
			&= u\bigg(\frac{d}{dz}P_{k}^{1}(u)\bigg)f^{(1)}(z) + \sum_{s=2}^{k}\bigg(uP_{k}^{s-1}(u) + u\frac{d}{dz}P_{k}^{s}(u)\bigg)f^{(s)}(z) + uP_{k}^{k}(u)f^{(k+1)}(z)
			.	\end{align*}
		
		Comparing the  left-hand side and the right-hand side of the last equality, 
		we obtain the 
		relations: 
		\begin{align}
			P_{k+1}^{1}(u) &= u\frac{d}{dz}P_{k}^{1}(u);\label{9}\\
			P_{k+1}^{s}(u) &= uP_{k}^{s-1}(u) + u\frac{d}{dz}P_{k}^{s}(u),
			\; 1 < s< k+1; \label{10}\\
			P_{k+1}^{k+1}(u) &= uP_{k}^{k}(u). \label{11}\\
			\intertext{Since $P_{1}^{1}(u)= u$, then from (\ref{11}) we obtain} 
			P_{k+1}^{k+1} &= u^{k+1} \text{\quad (i.e. (\ref{2}))}. \nonumber\\
			\intertext{From (\ref{9}) and (\ref{10}) we obtain the equalities: }
			P_{k+1}^{k} (u) &= uP_{k}^{k-1}(u) + ku^{k}\frac{d}{dz}u; \nonumber\\
			P_{k}^{k-1}(u) &= uP_{k-1}^{k-2}(u) + (k-1)u^{k-1}\frac{d}{dz}u;\nonumber \\
			&\;\;\vdots\nonumber\\
			P_{3}^{2}(u) & =  uP_{2}^{1}(u) + 2u^{2}\frac{d}{dz}u.\nonumber
		\end{align}
		We multiply the first equality by $1$, the second - by $u$, the third - by $u^{2}$, ..., the last - by $u^{k-2}$, then by summing the received equalities we obtain 
		\[
		P_{k+1}^{k}(u) = u^{k-1}P_{2}^{1}(u) + u^{k}\frac{d}{dz}u\sum_{t=2}^{k}t.
		\]
		
		Since $P_{2}^{1}(u) = u\frac{d}{dz}u$, we obtain 
		
		\begin{equation}\label{15-}
			P_{k+1}^{k}(u) = \bigg(\sum_{t=1}^{k}t\bigg)u^{k}\frac{d}{dz}u =  \frac{k(k+1)}{2}u^{k}\frac{d}{dz}u.
		\end{equation}
		
		To calculate $P_{k+1}^{k-1}(u)$ we use the same approach as for obtaining $(\ref{15-})$.
		
		Continuing the same way, one may also calculate $P_{k+1}^{k-2}(u)$, $P_{k+1}^{k-3}$, etc. and observe that in these particular cases (\ref{3}) is true.
		Below we shall prove (\ref{3}), using PMI.
		
		Substituting in (\ref{10}) $s$ with $k+1-s$, we obtain:
		\begin{equation}\label{13}
			P_{k+1}^{k+1-s}(u) = uP_{k}^{k-s}(u) + u\frac{d}{dz}P_{k}^{k+s+1}(u), \; s = 1, ..., k.
		\end{equation}
		
		With the help of (\ref{13}), we shall repeat the approach, realized for $P_{k+1}^{k}(u)$, but this time to obtain $P_{k+1}^{(k+1)-(s+1)}(u)$. Let us assume that $P_{k+1}^{k+1-s}(u)$ is given by (\ref{3}). Then
		
		\[
		\frac{d}{dz}P_{k+1}^{k+1-s}(u) = \sum_{m=1}^{s}(k+1-m)u^{k-m}F_{m}^{k+1}(u)\frac{d}{dz}u + \sum_{m=1}^{s}u^{k+1-m}\frac{d}{dz}F_{m}^{k+1}(u).
		\]
		
		From (\ref{13}), using the above equality, we obtain 
		\small
		\begin{align}\label{14}
			\frac{d}{dz}P_{k+1}^{k+1-s}(u) = \; &   u^{k}\frac{d}{dz}F_{1}^{k+1}(u) + \nonumber\\ & + \sum_{m=2}^{s}\bigg((k+2-m)F_{m-1}^{k+1}(u)\frac{d}{dz}u + \frac{d}{dz}F_{m}^{k+1}(u)\bigg)u^{k+1-m} +\nonumber\\ &  + (k+1-s)u^{k-s}F_{s}^{k+1}(u)\frac{d}{dz}u.
		\end{align}
		\normalsize	 
		Substituting in (\ref{14}) \, $k+1$ with $k, k-1, k-2... $, and after that multiplying the first equality by $u$, the second - by $u^{2}$, the third - by $u^{3}$, and so on, 
		then by summing the received equalities we obtain 
		
		\small
		\begin{align}\label{15}
			P_{k+1}^{(k+1)-(s+1)}(u) = \; &  \bigg(\sum_{t=s+1}^{k}\frac{d}{dz}F_{1}^{t}(u) \bigg)u^{k} \ + \nonumber\\
			& +\sum_{m=2}^{s}\bigg(\sum_{t=s+1}^{k}(t+1-m)F_{m-1}^{t}(u)\frac{d}{dz}u+ \frac{d}{dz}F_{m}^{t}(u) \bigg)u^{k+1-m} \ + \nonumber\\ 
			& + \bigg(\sum_{t=s+1}^{k}(t-s)F_{s}^{t}(u)\frac{d}{dz}u\bigg)u^{k-s}.
		\end{align}
		\normalsize
		Using the representation of $P_{k+1}^{(k+1)-(s+1)}(u)$, given in (\ref{15}), with the help of (\ref{3}) and (\ref{4}), we conclude that the relations (\ref{5}) - (\ref{8}) hold. Indeed, (\ref{4}) implies: 
		
		\begin{equation}\label{16}
			F_{m-1}^{t}(u)\frac{d}{dz}u = \sum_{m-1, \, s}  C_{\alpha_{1}, ..., \alpha_{s}}^{s, m-1, t}\bigg(\frac{du}{dz}\bigg)^{\alpha_{1} + 1 }\bigg(\frac{d^{2}u}{dz^{2}}\bigg)^{\alpha_{2} }...\bigg(\frac{d^{s}u}{dz^{s}}\bigg)^{\alpha_{s}}.
		\end{equation}
		On the other side, we have 
		\footnotesize
		\begin{align}\label{17}
			&\frac{d}{dz}F_{m}^{t}(u) = \nonumber\\ \;\;&\sum_{m,s}\sum_{j=1}^{s}\delta_{j}C_{\delta_{1}, ..., \delta_{s}}^{s,m,t}\bigg(\frac{du}{dz}\bigg)^{\delta_{1}}...\bigg(\frac{d^{j-1}u}{dz^{j-1}}\bigg)^{\delta_{j-1}}\bigg(\frac{d^{j}u}{dz^{j}}\bigg)^{\delta_{j}-1}\bigg(\frac{d^{j+1}u}{dz^{j+1}}\bigg)^{\delta_{j+1}+1}\bigg(\frac{d^{j+2}u}{dz^{j+2}}\bigg)^{\delta_{j+2}}... \bigg(\frac{d^{s}u}{dz^{s}}\bigg)^{\delta_{s}}.
		\end{align}
		\normalsize
		
		Comparing (\ref{1}), (\ref{3}), (\ref{4}) with (\ref{15}) - (\ref{17}), we obtain: 
		\[
		\delta_{1} = \alpha_{1} + 1, \dots , \ \delta_{j}-1= \alpha_{j}, \  \delta_{j+1}+1=\alpha_{j+1}, \dots, \ \delta_{s}=\alpha_{s}.
		\]
		Thus, we get (\ref{5}). The relations (\ref{6}) and (\ref{7}) could be received analogically. To finish the proof of Theorem 3, it remains to prove that each monomial in $\displaystyle\sum_{m,s+1}$, for example
		
		\[\bigg(\frac{du}{dz}\bigg)^{\gamma_{1} }\bigg(\frac{d^{2}u}{dz^{2}}\bigg)^{\gamma_{2} }...\bigg(\frac{d^{s}u}{dz^{s}}\bigg)^{\gamma_{s}}\bigg(\frac{d^{s}u}{dz^{s+1}}\bigg)^{\gamma_{s+1}}\]
		with $\gamma_{s+1}=0$, is contained in the right-hand side of (\ref{16}). This fact can be proved in the same way as in the proof of Lemma 3.2 from \cite{Petkova}. The theorem is proved.
	\end{proof}
	
	\bigskip
	
	Another approach for finding the expansion of $A^{k}$ is proposed below. 
	
	We set: $A^{0} = u(z)$; $u(z) = \frac{1}{v(z)}$. Thus $A = \frac{d}{d\int v(z)}$. Let $\int v(z) = \varphi(z)$. Then $A = \frac{d}{d\varphi(z)}$. Therefore, $Af(z) = \frac{d}{d\varphi(z)}f(z)$. Substituting $\varphi(z) = t$, we have $z = \varphi ^{-1}(t)$ (where  $\varphi ^{-1}$ is the inverse function of $\varphi$). Therefore
	\[
	Af(z) = Af(\varphi^{-1}(t)) = \frac{d}{dt}f(\varphi^{-1}(t))|_{t=\varphi(z)}.
	\]
	Hence
	\[
	A^{k}f(z) = \left(\frac{d}{dt}\right)^{k}f(\varphi^{-1}(t)) | _{t = \varphi(z)}, \; k=1, 2, 3,  \dots
	\]
	Now, for the right-hand side of the last equality Fa
	\`{a} di Bruno's formula (\cite{Abr}, p.823) is applicable and as a result we obtain
	\begin{equation}\label{18}
		A^{k}f(z) = \sum_{m=0}^{k} f^{(m)}(\varphi^{-1}(t))\sum_{m, k} C_{k}(\alpha)\left[[D  \varphi^{-1}(t)]^{\alpha_{1}}\dots [D^{k} \varphi^{-1}(t)]^{\alpha_{k}}\right]|_{t = \varphi(z)},
	\end{equation}
	where $\alpha = (\alpha_{1}, \alpha_{2}, \dots, \alpha_{k})$ and \begin{equation}\label{19}
		C_{k}(\alpha) = \frac{k!}{(1!)^{\alpha_{1}}\dots (k!)^{\alpha_{k}}\alpha_{1}!\dots \alpha_{k}!}.
	\end{equation}
	
	From the formula for derivatives of inverse functions we obtain:
	\begin{align*}
		(\varphi^{-1}(t))' &= \frac{1}{\varphi'(z)} = u(z) = A^{0}u; \\
		(\varphi^{-1}(t))'' &= u'(z) (\varphi^{-1}(t))' = u(z)u'(z) = A^{1}u; \dots;\\
		(\varphi^{-1}(t))^{(k)} &= A^{k-1}u.
	\end{align*}
	
	Substituting in (\ref{18}) the results of the above equalities, we obtain:
	\begin{equation}\label{20}
		A^{k}f(z) = \sum_{m=0}^{k} f^{(m)}(z)\sum_{m, k} C_{k}(\alpha)[A^{0}u]^{\alpha_{1}}\dots [A^{k-1}u]^{\alpha_{k}}. 
	\end{equation}
	Thus (using (\ref{20})), we proved the following result:
	
	\begin{theorem}\label{thm:theorem2.2}
		The operator $A^{k}$, $k=1,2, \dots$, could be given by:
		\begin{equation}\label{21}
			A^{k} = \sum_{m=0}^{k} \sum_{m, k} C_{k}(\alpha)[A^{0}u]^{\alpha_{1}}\dots [A^{k-1}u]^{\alpha_{k}}\left(\frac{d}{dz}\right)^{m}.
		\end{equation}
	\end{theorem}

	\section{Applications of the obtained results}
	
	\begin{corollary}\label{cor:corollary3.1}
		It is valid
		\begin{equation}\label{22}
			C_{\underbrace{\text{\scriptsize $0, ..., 0$}}_{\text{\normalfont s-1}}, 1}
			^{s, 1, k+1} = {k+1 \choose s+1}, s = 1, ..., k.
		\end{equation}
	\end{corollary}
	\begin{proof}
		For $s=1$ we have $\displaystyle C_{1}
		^{1, 1, k+1} = {k+1 \choose 2}.$ Then, to prove (\ref{22}), we must verify that
		
		\[\sum_{t=s+1}^{k}{t \choose s+1} = {k+1\choose s+1}.\]   
		But the last identity follows immediately from the definition of the binomial coefficients: \[{p\choose q} = \frac{p!}{q!(p-q)!} \cdot \]
	\end{proof}
	
	\begin{corollary}\label{cor:corollary3.2}
		It is valid
		
		\begin{equation}\label{23}
			C_{s+1,\;\underbrace{\text{\scriptsize $0, ..., 0$}}_{\text{\normalfont s}}}
			^{s+1, s+1, k+1} = \sigma(k+1, k-s).
		\end{equation}
	\end{corollary}
	\begin{proof} For $s=0$  (\ref{23}) is true. Indeed, the left-hand side of  (\ref{23}) for $s=0$ is: 
		\[C_{1}^{1, 1, k+1} = {k+1 \choose 2} \quad (\text{from } (\ref{22}) \text{ with } s=1).
		\]
		But \ $\sigma(k+1,k) = {k+1 \choose 2}$ \ (see \cite{Abr}, p.~825). Thus, 
		\[
		C_{1}^{1, 1, k+1} = \sigma(k+1, k).
		\]
		
		Further, we shall prove (\ref{23}), using PMI. Let us denote
		\begin{equation}\label{24}
			C_{s+1,\;\underbrace{\text{\scriptsize $0, ..., 0$}}_{\text{\normalfont s}}}
			^{s+1, s+1, k+1} = a(k+1, k-s), \quad s=0, 1, 2, \dots, k; \ k\geq0.
		\end{equation} 
		Then, from (\ref{6}),
		\begin{equation*}
			a(k+1,k-s) = \sum_{t=s+1}^{k}(t-s)a(t,t-s).
		\end{equation*}
		Let $k-s=\tilde{s}$. Then the above equality yields
		\begin{equation}\label{25}
			a(k+1, \tilde{s}) = \sum_{t=k+1-\tilde{s}}^{k}\left(t+\tilde{s}-k\right)a(t,t+\tilde{s}-k),
		\end{equation}
		$\tilde{s} = 0, 1, \dots, k, k \geq 0$.
		
		To prove that $a(k, \tilde{s}) = \sigma(k,\tilde{s})$, it is enough to establish that: 
		
		\begin{align}
			a(k+1,\tilde{s}) &= a(k,\tilde{s}-1) + \tilde{s}a(k, \tilde{s})\label{26}\\
			a(k+1,k) & = {k+1 \choose 2}, \label{27}
		\end{align}
		since the same recurrence relations are satisfied by the Stirling numbers of the second kind (see \cite{Abr}, p.~825).
		
		From the definition of $a(k+1, k-s)$ we obtain that $a(k+1, \tilde{s}) = C_{k+1-\tilde{s},\;\underbrace{\text{\scriptsize $0, ..., 0$}}_{\text{\normalfont $k-\tilde{s}$}}}^{k+1-\tilde{s},\, k+1-\tilde{s}, \, k+1}$.
		
		For $\tilde{s}=k$, the above equality yields $a(k+1, k) = C_{1}^{1, \,1,\, k+1} = {k+1 \choose 2}$ and (\ref{27}) is proved.
		
		With the help of (\ref{25}), the validity of (\ref{26}) is a matter of direct check. 
		
		Thus, we not only proved (\ref{23}), but also obtained the following recurrence relation for the Stirling numbers of the second kind:
		\begin{equation}\label{28}
			\sigma(k+1, k-s) = \sum_{t=s+1}^{k} (t-s) \sigma(t, t-s), 
		\end{equation}
		$s=0, \dots, k; \ \sigma(k, k-1) = {k \choose 2}$. 
		
		One may use (\ref{28}) as another defintion for the Stirling numbers of the second kind.
		
		Since $\sigma(k,m)$ has the explicit representation (see \cite{Abr}, p.~824):
		
		\[
		\sigma(k,m) = \frac{1}{m!}\sum_{i=0}^{m}(-1)^{m-i}{m \choose i}i^{k},
		\]
		then we obtain the numbers $C_{s+1,\;\underbrace{\text{\scriptsize $0, ..., 0$}}_{\text{\normalfont s}}}
		^{s+1, s+1, k+1}$ in explicit form:
		\[
		C_{s+1,\;\underbrace{\text{\scriptsize $0, ..., 0$}}_{\text{\normalfont s}}}
		^{s+1, s+1, k+1} = \frac{1}{(k-s)!}\sum_{t=0}^{k-s}(-1)^{k-s-t}{k-s \choose t}t^{k+1}.
		\]
	\end{proof}	
	
	From (\ref{1}), if we change $k+1$ with $k$ and define that $P_{k}^{0}(u) = 0$, we obtain 
	\begin{equation}\label{ak}
		A^{k} = \sum_{s=0}^{k}P_{k}^{s}(u)\bigg(\frac{d}{dz}\bigg)^{s}.
	\end{equation} If we put $k-s$ instead of $s$ in (\ref{ak}), we obtain
	\begin{equation} 
		\label{32}
		A^{k} = \sum_{s=0}^{k}P_{k}^{k-s}(u)\bigg(\frac{d}{dz}\bigg)^{k-s}.
	\end{equation}
	
	Now we replace in (\ref{32}) $P_{k}^{k-s}(u)$ with the right-hand side of (\ref{3}) with $k$ instead $k+1$. As a result, we obtain:
	
	\begin{equation}\label{33}
		A^{k} = \sum_{s=0}^{k}\left(\sum_{m=1}^{s}u^{k-m}\sum_{m,s} C_{\alpha_{1}, ..., \alpha_{s}}^{s, m, k}\bigg(\frac{du}{dz}\bigg)^{\alpha_{1}}...\bigg(\frac{d^{s}u}{dz^{s}}\bigg)^{\alpha_{s}} \right) \left(\frac{d}{dz}\right)^{k-s}.
	\end{equation}
	
	Using that $\alpha_{1}+\dots+\alpha_{s} = m$, one may rewrite (\ref{33}) in the form:
	\begin{equation}\label{33''}
		A^{k} = u^{k}\sum_{s=0}^{k}\left(\sum_{m=1}^{s}\sum_{m,s} C_{\alpha_{1}, ..., \alpha_{s}}^{s, m, k}\bigg(\frac{u'}{u}\bigg)^{\alpha_{1}}...\bigg(\frac{u^{(s)}}{u}\bigg)^{\alpha_{s}} \right) \left(\frac{d}{dz}\right)^{k-s}.
	\end{equation}
	
	The representation (\ref{33''}) gives 
	the possibility to express the powers of the operator $B=e^{g(z)}\frac{d}{dz}$ with the help of the complete $n\!-\!th$ Bell polynomial $B_{n}(g_{1}, \dots, g_{n})$, where $g_{i} = D^{i}g$, $i=1,\dots,n$, and $D=\frac{d}{dz}$. Namely, setting in $(35)$ $u(z) = e^{g(z)}$, we obtain 
	\begin{equation}\label{33'''}
		B^{k} = \left(e^{g(z)}\right)^{k}\sum_{s=0}^{k}\left(\sum_{m=1}^{s}\sum_{m,s} C_{\alpha_{1}, ..., \alpha_{s}}^{s, m, k}
		\prod_{n=1}^{s}\left(e^{-g(z)}D^{n}e^{g(z)}\right)^{\alpha_{n}}\right)
		D^{k-s}.
	\end{equation}
	
	Thus, we proved 
	
	\begin{proposition}	\label{prop:proposition3.1}
		The operator $B^{k}$, $k\geq1$, admits the representation
		\begin{equation}\label{33''''}
			B^{k} = \left(e^{g(z)}\right)^{k}\sum_{s=0}^{k}\left(\sum_{m=1}^{s}\sum_{m,s} C_{\alpha_{1}, ..., \alpha_{s}}^{s, m, k}
			\prod_{n=1}^{s}\left(B_{n}(g_{1}, \dots, g_{n})\right)^{\alpha_{n}}\right)
			\left(\frac{d}{dz}\right)^{k-s}.
		\end{equation}
		
		Indeed, we have $B_{n}(g_{1}, \dots, g_{n}) =  e^{-g(z)}D^{n}e^{g(z)}$ (see \cite{wiBell}). \hfill
		$	{\square}$
	\end{proposition}
	\bigskip
	
	Now, we shall use (\ref{33}) to obtain some corollaries.
	
	\begin{corollary}\label{cor:corollary3.3}
		The identity 
		
		\begin{equation}
			\sum_{m=1}^{s}\sum_{m,s}C_{\alpha_{1}, ..., \alpha_{s}}^{s,m,k} = (-1)^{s}S(k,k-s), \quad s=0, \dots, k, \label{34}
		\end{equation}	
		holds.
	\end{corollary}
	
	\begin{proof}
		First we put in (\ref{b}) $k-s$ instead of $s$ to obtain:
		\begin{equation} \label{35}
			A_{2} \! ^{k} =  e^{kz}\sum_{s=0}^{k}(-1)^{s}S(k,k-s)\left(\frac{d}{dz}\right)^{k-s}.
		\end{equation}
		
		Second, we put $u(z) = e^{z}$ in (\ref{33}) and using that $\alpha_{1}+\alpha_{2}+\dots+\alpha_{s} = m$, we  obtain 
		
		\begin{equation}\label{36}
			A_{2} \! ^{k} =  e^{kz}\sum_{s=0}^{k}\left(\sum_{m=1}^{s}\sum_{m,s} C_{\alpha_{1}, ..., \alpha_{s}}^{s, m, k} \right)\left(\frac{d}{dz}\right)^{k-s}. 
		\end{equation}
		
		Now we compare (\ref{35}) and (\ref{36}) and thus (\ref{34}) is proved.
	\end{proof}
	
	\bigskip
	
	\begin{corollary}\label{cor:corollary3.4}
		The identity 
		\begin{equation}\label{37}
			\sum_{m=1}^{s}\sum_{m,s}(1!)^{\alpha_{1}}...(s!)^{\alpha_{s}}C_{\alpha_{1}, ..., \alpha_{s}}^{s,m,k} = (2s-1)!!{k+s-1 \choose 2s}
		\end{equation}
		holds. 
	\end{corollary}
	\begin{proof}
		In (\ref{c}) we put $k-s$ instead of $s$ and obtain
		
		\begin{equation}\label{38}
			A_{3}\!^{k} = \sum_{s=0}^{k}(-1)^{s}(2s-1)!!{k+s-1 \choose 2s}\frac{1}{z^{k+s}}\left(\frac{d}{dz}\right)^{k-s}.
		\end{equation}
		Let us put $u(z) = \frac{1}{z}$ in (\ref{33}). Then we obtain
		
		\begin{align*}
			A_{3}\!^{k} &= \sum_{s=0}^{k}\left(\sum_{m=1}^{s}\frac{1}{z^{k-m}}\sum_{m,s} C_{\alpha_{1}, ..., \alpha_{s}}^{s, m, k}\bigg(\frac{dz^{-1}}{dz}\bigg)^{\alpha_{1}}...\bigg(\frac{d^{s}z^{-1}}{dz^{s}}\bigg)^{\alpha_{s}} \right) \left(\frac{d}{dz}\right)^{k-s} = \\
			&= \sum_{s=0}^{k}\left(\sum_{m=1}^{s}\frac{1}{z^{k-m}}\sum_{m,s} C_{\alpha_{1}, ..., \alpha_{s}}^{s, m, k}\bigg(\frac{1!(-1)^{1}}{z^{2}}\bigg)^{\alpha_{1}}...\bigg(\frac{s!(-1)^{s}}{z^{s+1}}\bigg)^{\alpha_{s}} \right) \left(\frac{d}{dz}\right)^{k-s} = \\
			&= \sum_{s=0}^{k}\left(\sum_{m=1}^{s}\frac{1}{z^{k-m}}\sum_{m,s} (-1)^{\sum_{i=1}^{s}i\alpha_{i}}\prod_{i=1}^{s}(i!)^{\alpha_{i}} C_{\alpha_{1}, ..., \alpha_{s}}^{s, m, k}\frac{1}{z^{\sum_{i=1}^{s}(i+1)\alpha_{i}}}
			\right) \left(\frac{d}{dz}\right)^{k-s}.
		\end{align*}
		Since: $\displaystyle\sum_{i=1}^{s}i\alpha_{i} = s$; $\displaystyle\sum_{i=1}^{s}\alpha_{i} = m$, we have:
		\[
		\sum_{i=1}^{s}(i+1)\alpha_{i} = \sum_{i=1}^{s}\alpha_{i} + \sum_{i=1}^{s}i\alpha_{i} = m+s.
		\]
		Therefore, 
		\begin{equation}\label{39}
			A_{3}\!^{k} 
			= \sum_{s=0}^{k}(-1)^{s}\left(\sum_{m=1}^{s}\sum_{m,s} \prod_{i=1}^{s}(i!)^{\alpha_{i}} C_{\alpha_{1}, ..., \alpha_{s}}^{s, m, k}\right)\frac{1}{z^{k+s}}
			\left(\frac{d}{dz}\right)^{k-s}.
		\end{equation}
		Comparing the last equality with (\ref{38}), we obtain (\ref{37}).
	\end{proof}
	\bigskip
	
	\begin{corollary}\label{cor:corollary3.5}
		The identity 
		\begin{equation}\label{40}
			\sum_{s,k}\frac{k!}{1^{\alpha_{1}}\alpha_{1}!...k^{\alpha_{k}}.\alpha_{k}!} = \sum_{m=1}^{k-s}\sum_{m,k-s}C_{\alpha_{1}, ..., \alpha_{k-s}}^{k-s,m, k}
		\end{equation}
		holds. 
	\end{corollary}
	\begin{proof}
		The relation 
		\begin{equation}\label{41}
			\sum_{s,k}\frac{k!}{1^{\alpha_{1}}\alpha_{1}!...k^{\alpha_{k}}.\alpha_{k}!} = (-1)^{k-s}S(k,s)
		\end{equation}
		is well-known (see (\cite{Abr}, p.~823)). 
		From Corollary 3 (see (\ref{34})):
		\begin{equation}
			\sum_{m=1}^{s}\sum_{m,s}C_{\alpha_{1}, ..., \alpha_{s}}^{s,m,k} = (-1)^{s}S(k,k-s), \quad s=0, \dots, k. \label{42}
		\end{equation}	
		In (\ref{42}) we put $k-s$ instead of $s$ to obtain:
		\begin{equation}
			\sum_{m=1}^{k-s}\sum_{m,k-s}C_{\alpha_{1}, ..., \alpha_{k-s}}^{k-s,m,k} = (-1)^{k-s}S(k,s), \quad s=0, \dots, k \label{43}.
		\end{equation}	
		Comparing the last equality with (\ref{41}), we prove (\ref{40}).
	\end{proof}
	
	\section{A correspondence between the operators \boldmath$A_{1} = A_{1, x} = x\frac{d}{dx}$ and $A= A_{z} = u(z)\frac{d}{dz}$}\unboldmath
	
	Let the substitution $x=\varphi(z)$ transform $A_{1,x}$ into $A_{z}$. Then \begin{equation} \label{44}
		\frac{\varphi(z)}{\varphi'(z)} = u(z).
	\end{equation}
	The equation (\ref{44}) has the solution $\varphi(z) = e^{\int_{c}^{z}\frac{dt}{u(t)}}$, where $c$ is an arbitrary complex constant. Thus, $\varphi$ transforms $\displaystyle A_{1,x}\!^{k} = \sum_{s=0}^{k}\sigma(k,s)x^{s}\left(\frac{d}{dx}\right)^{s}$ into $\displaystyle A_{z}\!^{k} = \sum_{s=0}^{k}P_{k}^{s}(u)\left(\frac{d}{dz}\right)^{s}$, $k=1, 2, \dots$ (see (\ref{a}) and (\ref{ak})). Now, for an arbitrary function $f$, on which the operator $A_{z}\!^{s}$ acts,  we apply Fa
	\`{a} di Bruno's formula:
	\begin{equation}\label{45}
		\left(\frac{d}{dx}\right)^{s}f(z(x)) = \sum_{m=0}^{s}f^{(m)}(z(x))\Delta(s,m), \; s=1,2,3, \dots, 
	\end{equation}
	where
	
	\[\Delta(s,m) = \sum_{m,s}C_{s}(\alpha)\left(\frac{d}{dx}z(x)\right)^{\alpha_{1}} \dots \left( \left(\frac{d}{dx}\right)^{s}z(x)\right)^{\alpha_{s}} \; \text{(see (\ref{19}))}.
	\]
	Thus, the operator $D^{s}$ transforms into the operator $\displaystyle\sum_{m=1}^{s}\Delta(s,m)\left(\frac{d}{dz}\right)^{m}$ and moreover we have the identities:
	
	\begin{equation}\label{47''}
		\sum_{i=0}^{k-s}\sigma(k,s+i)\Delta(s+i,s) (\varphi(z))^{s+i} 
		= P_{k}^{s}(u), 
	\end{equation}
	$s=1, \dots, k$. 
	
	One may rewrite the above identities in the form:
	\begin{equation}\label{47'}
		\sigma(k,s)V(s,s) x^{s} + \sigma(k,s+1)V(s+1,s)x^{s+1} + \dots +	\sigma(k,k)V(k,s) x^{k} = P_{k}^{s}\left(u\left(\varphi^{-1}(x)\right)\right),
	\end{equation}
	$s=1, \dots, k-1$; $P_{k}^{k}(u) = u^{k}\left(\varphi^{-1}(x)\right)$, where $\varphi^{-1}$ is the inverse function of $\varphi$ and 
	\[
	V(n,s) = \sum_{s,n}C_{n}(\alpha)\left(\frac{d}{dx} \ \varphi^{-1}(x)\right)^{\alpha_{1}} \dots \left( \left(\frac{d}{dx}\right)^{n} \varphi^{-1}(x)\right)^{\alpha_{n}}. \;
	\]
	
	Let $s=1$ in (\ref{47'}). Since $\alpha_{1} + \alpha_{2} + \dots +\alpha_{n} = s = 1$, then $ \alpha_{1} = \alpha_{2} = \dots =\alpha_{n-1} = 0$;   $\alpha_{n} = 1$; $C_{n}(\alpha)=1$. Hence
	\[
	V(n,s) = V(n,1) = \left(\frac{d}{dx}\right)^{n}\varphi^{-1}(x).
	\]
	
	Therefore, (\ref{47'}) yields:
	
	\begin{corollary}
		\label{cor:corollary4.1}
		The identity
		\begin{equation}\label{48}
			\sum_{i=1}^{k}\sigma(k,i)\left(\left(\frac{d}{dx}\right)^{i}\varphi^{-1}(x)\right) x ^{i} = P_{k}^{1}\left( u\left(\varphi^{-1}(x)\right)\right)
		\end{equation}
		holds, where $k$ is an arbitrary natural number. \hfill${\square}$
	\end{corollary}
	
	It is obvious that the function $\psi(x) = \varphi^{-1}(x)$ satisfies the differential equation 	 
	\begin{equation}\label{diff}
		x\psi'(x) = u(\psi(x)).
	\end{equation}

	For the case $u(z) = \frac{1}{z}$, (\ref{diff}) yields $\psi(x) = \sqrt{\ln(x^{2})}$. From (\ref{48}) and (\ref{c}) (for $s=1$) it follows:
	
	\begin{corollary}\label{4.2}
		For an arbitrary natural number $k$ the identity
		\footnotesize
		\begin{align}
			&\sigma(k,1) \left(\frac{d}{dx}\sqrt{ln(x^{2})}\right)x +  \sigma(k,2) \left(\left(\frac{d}{dx}\right)^{2}\sqrt{ln(x^{2})}\right)x^{2} + \dots +  \sigma(k,k) \left(\left(\frac{d}{dx}\right)^{k}\sqrt{ln(x^{2})}\right)x^{k} =\nonumber\\[5pt]
			&=(-1)^{k-1}(2k-3)!!\left(\sqrt{\ln(x^{2})}\right)^{1-2k}, \nonumber
		\end{align}
		\normalsize
		holds, where we set $(-1)!!=1$.
		\hfill$\square$
	\end{corollary}
	
	The following lemma is the key for finding $\left(\frac{d}{dx}\right)^{n}\psi(x)$, $n=1, 2, 3, \dots$, where $\psi$ is 
	the solution of (\ref{diff}).
	\bigskip
	
	\begin{lemma}\label{lem:lemma4.1}
		Let $T_{1}(z) = u(z)$ and $T_{n+1}(z) = u(z)\frac{d}{dz}T_{n}(z) - nT_{n}(z)$, $n\geq1$.  
		Then 
		\[
		\left(\frac{d}{dx}\right)^{n}z = \left(\frac{d}{dx}\right)^{n}\psi(x) = \frac{ T_{n}(z)}{(\varphi(z))^{n}}.
		\]
	\end{lemma}
	\begin{proof}
		Since $z = \varphi^{-1}(x) = \psi(x)$, then
		\[
		\frac{d}{dx}z = \frac{d}{dx}\varphi^{-1}(x) = \frac{1}{\varphi'(z)}  \overset{\text{from (\ref{44})}}{=} \frac{u(z)}{\varphi(z)} =  \frac{T_{1}(z)}{\varphi(z)},
		\]
		which proves the lemma for $n=1$. 
		
		Let the assertion of the lemma be true for some $n\geq1$. Then
		\begin{align*}
			\left(\frac{d}{dx}\right)^{n+1}z \ &= \ \frac{d}{dx} \left(\frac{d}{dx}\right)^{n} z \ \overset{\text{from PMI}}{=} \frac{d}{dx} \frac{T_{n}(z)}{(\varphi(z))^{n}} =\\[5pt]
			&= \frac{\left(\varphi(z)\right)^{n}\left(\frac{d}{dz}T_{n}(z)\right)\frac{d}{dx}z - nT_{n}(z) \left(\varphi(z)\right)^{n-1}\frac{d}{dx}\varphi(z)}{\left(\varphi(z)\right)^{2n}} =\\
			&= \frac{\left(\varphi(z)\right)^{n-1}\left[u(z)\frac{d}{dz}T_{n}(z
				) - nT_{n}(z)\left(\frac{d}{dz}\varphi(z)\right)\frac{d}{dx}z\right]}{\left(\varphi(z)\right)^{2n}} =\\
			&= \frac{T_{n+1}(z)}{\left(\varphi(z)\right)^{n+1}}
		\end{align*}
		and the lemma is proved. 
	\end{proof}
	
	\begin{remark}\label{rem:remark4.1}
		It is obvious that  
		\[
		\left(\frac{d}{dz}\varphi(z)\right)\frac{d}{dx}z =1,  
		\]
		since \ $\displaystyle\varphi'(z) = \frac{\varphi(z)}{u(z)}$ \ and \ $\displaystyle\frac{d}{dx}z = \frac{u(z)}{\varphi(z)}$.
		\; \hfill${\square}$
	\end{remark}
	\bigskip

	Let $u(z) = e^{z}$. Then   
	\begin{equation}\label{50'}
		T_{n}(z) = r(n,n)\left(e^{z}\right)^{n} + r(n,n-1)\left(e^{z}\right)^{n-1} + \dots + r(n,1)\left(e^{z}\right)^{1},
	\end{equation}
	where the numbers $r(n,k)$ satisfy the recurrence relations:
	\begin{align}\label{52'}
		r(n+1, n+1) &= nr(n,n) = n!; \nonumber\\
		r(n+1,1) &= -nr(n,1) = (-1)^{n}n!;\\
		r(n+1,k) &= (k-1)r(n,k-1) - nr(n,k), \quad k=2, \dots, n. \nonumber
	\end{align}
	Now we use the lemma (\! for $u(z) = e^{z}$\! ) to obtain:
	\[
	\Delta(s,m) = \sum_{m,s} C_{s}(\alpha) \prod_{n=1}^{s} \left[\frac{T_{n}(z)}{(\varphi(z))^{n}}\right]^{\alpha_{n}} =  (\varphi(z))^{-(\alpha_{1} +2\alpha_{2} + \dots+ n\alpha_{n})} \sum_{m,s} C_{s}(\alpha) \prod_{n=1}^{s} \left[T_{n}(z)\right]^{\alpha_{n}}. 
	\]
	Hence, from (\ref{50'}):
	\begin{equation}\label{53}
		\Delta(s,m) = \left(\varphi(z)\right)^{-s}\sum_{m,s}C_{s}(\alpha)\prod_{n=1}^{s} \left[\sum_{i=1}^{n}r(n,i)\left(e^{z}\right)^{i}\right]^{\alpha_{n}},
	\end{equation}
	because $\alpha_{1} + 2\alpha_{2} + \dots + s\alpha_{s} = s $.
	
	From (\ref{47''}) and (\ref{53}) we obtain:
	\begin{equation*}
		P_{k}^{s}(e^{z}) = \sum_{i=0}^{k-s} \sigma(k,s+i)\sum_{s,s+i}C_{s+i}(\alpha)\prod_{n=1}^{s+i}\left(\sum_{j=1}^{n}r(n,j)\left(e^{z}
		\right)^{j}\right)^{\alpha_{n}}.
	\end{equation*}
	
	On the other hand, from (\ref{b}), we have:
	\[
	P_{k}^{s} \left(e^{z}\right) = (-1)^{k+s}S(k,s)e^{kz},
	\]
	where 
	$S(k,s)$, as before, are the Stirling numbers of the first kind.
	
	The last two equalities yield:
	\begin{equation}\label{pks}
		(-1)^{k+s}S(k,s)e^{kz} = \sum_{i=0}^{k-s} \sigma(k,s+i)\sum_{s,s+i}C_{s+i}(\alpha)\prod_{n=1}^{s+i}\left(\sum_{j=1}^{n}r(n,j)\left(e^{z}
		\right)^{j}\right)^{\alpha_{n}}.
	\end{equation}
	We set $z=0$ in (\ref{pks}) to obtain
	\begin{equation}\label{pks2}
		(-1)^{k+s}S(k,s) = \sum_{i=0}^{k-s} \sigma(k,s+i)\sum_{s,s+i}C_{s+i}(\alpha)\prod_{n=1}^{s+i}R_{n}\!^{\alpha_{n}},
	\end{equation}
	where $\displaystyle R_{n} = \sum_{j=1}^{n} r(n,j) = T_{n}(0)$.
	
	Let us remind that $(-1)^{k+1}S(k,1) = (k-1)!$ (see \cite{Abr}, p.~824). Therefore, setting $s=1$  in (\ref{pks}), we obtain:
	\begin{equation}\label{56'}
		(k-1)!e^{kz} = \sum_{i=0}^{k-1}\sigma(k,i+1)\sum_{1, i+1}C_{i+1}(\alpha)\prod_{n=1}^{i+1}\left(\sum_{j=1}^{n}r(n,j)\left(e^{z}\right)^{j}\right)^{\alpha_{n}} = P_{k}^{1}\left(e^{z}\right).
	\end{equation}
	Since $x^{i} = \left(\varphi(z)\right)^{i}$, $i=1, 2, \dots, k$, then the lemma and (\ref{48}) yield:
	
	\begin{equation}\label{57'}
		P_{k}^{1}\left(e^{z}\right) = \sum_{i=1}^{k}\sigma(k,i)T_{i}(z).
	\end{equation}
	From (\ref{50'}), (\ref{56'}) and (\ref{57'}) we obtain 
	
	\[
	(k-1)!e^{kz} = \sum_{i=1}^{k}\sigma(k,i)\sum_{\nu=1}^{i}r(i,\nu)\left(e^{z}\right)^{\nu}.
	\]
	Since $\sigma(k,k) = 1$ and $r(k,k) = (k-1)!$ (see (\ref{52'})), then for $k\geq2$ one may rewrite the above equality in the form:
	\begin{equation}\label{58'}
		\sum_{n=1}^{k-1}\left[\sigma(k,n)r(n,n) + \sigma(k,n+1)r(n+1,n) + \dots + \sigma(k,k)r(k,n)\right]\left(e^{z}\right)^{n} = 0.
	\end{equation}
	
	From (\ref{58'}) we obtain:
	
	\begin{corollary}\label{cor:corollary4.3}
		For $k\geq2$ the relations for ortogonality:
		\begin{equation}\label{59'}
			\sigma(k,n)r(n,n) + \sigma(k,n+1)r(n+1,n) + \dots + \sigma(k,k)r(k,n) = 0,
		\end{equation}
		$n=1, \dots, k-1$, hold. \hfill $\square$
	\end{corollary}
	
	Let $r(n,k) = d_{k}S(n,k)$ (where again $S(n,k)$ are the Stirling numbers of the first kind). Replacing this equality in (\ref{52'}) we obtain that $d_{1} = 1$ and $d_{k}  =(k-1)d_{k-1}$, $k\geq2$. Therefore $d_{k}  = (k-1)!$. Hence, 
	\begin{equation}\label{60'}
		r(n,k) = (k-1)!S(n,k), \quad k=1, \dots, n.
	\end{equation} 
	For $k\geq2$, from 	(\ref{59'}) and (\ref{60'})  the well-known relations for ortogonality between the Stirling numbers of the first and of the second kind hold immediately:
	
	\begin{equation}\label{61'}
		\sigma(k,n)S(n,n) + \sigma(k,n+1)S(n+1,n) + \dots + \sigma(k,k)S(k,n) = 0,
	\end{equation}
	$n=1, \dots, k-1$.
	
	Using (\ref{60'}), we obtain $\displaystyle R_{n} = \sum_{\nu=1}^{n}(\nu-1)!S(n,\nu)$ and then (\ref{pks2}) yields{ 
		
		\begin{corollary}\label{cor:corollary4.4}
			The identity
			\[
			(-1)^{k+s}S(k,s) = \sum_{i=0}^{k-s}\sigma(k,s+i)\sum_{s, s+i}
			C_{s+i}(\alpha)\prod_{n=1}^{s+i}\left(\sum_{j=1}^{n}(j-1)!S(n,j)\right)^{\alpha_{n}}
			\]
			holds.	\hfill	$\square$
		\end{corollary}
		
		Let us consider the operator $A_{3} = \frac{1}{z}\frac{d}{dz}$. Obviously, $A_{3} = 2\frac{d}{dz^{2}}$. Hence,
		\[
		A_{3}\!^{k+1} = 2^{k+1}\left(\frac{d}{dz^{2}}\right)^{k+1}.
		\]
		
		We set $z=\sqrt{y}$ \ to obtain  
		\[
		\left(A_{3}\!^{k+1}(\sqrt{y})\right)f = 2^{k+1}f^{(k+1)}(\sqrt{y}).
		\]
		
		Further, we apply Fa
		\`{a} di Bruno's formula for calculating $f^{(k+1)}(\sqrt{y})$ to obtain:
		
		\[
		\left(A_{3}\!^{k+1}(\sqrt{y})\right)f = 2^{k+1}\sum_{m=0}^{k+1}f^{(m)}(\sqrt{y})\sum_{m,k+1}C_{k+1}(\alpha)\prod_{s=1}^{k+1}\left((\sqrt{y})^{(s)}\right)^{\alpha_{s}}.
		\]
		After further calculation of the product  $\displaystyle\prod_{s=1}^{k+1}\left((\sqrt{y})^{(s)}\right)^{\alpha_{s}}$ and after changing $\sqrt{y}$ with $z$, we finally obtain:
		
		\begin{equation} \label{62'}
			A_{3}\!^{k+1} = 2^{k+1}\sum_{s=0}^{k}\frac{1}{z^{k+1-s}}\left(\sum_{s,k+1}C_{k+1}(\alpha)\prod_{\nu=1}^{k+1}(\nu!)^{\alpha_{\nu}}{ \frac{1}{2}\choose \nu}^{\alpha_{\nu}}\right)\left(\frac{d}{dz}\right)^{k+1-s}.
		\end{equation}
		
		After replacing $k+1$ with $k$ in (\ref{38}) and comparing with (\ref{62'}), we obtain: 
		
		\begin{corollary}\label{cor:corollary4.5}
			The identities:
			\begin{equation*}
				2^{k+1}	\sum_{s,k+1}C_{k+1}(\alpha)\prod_{\nu=1}^{k+1}(\nu!)^{\alpha_{\nu}}{ \frac{1}{2}\choose \nu}^{\alpha_{\nu}} = (-1)^{s}(2s-1)!!{k+s\choose 2s},
			\end{equation*}
			$s=1,\dots,k$, hold, where we set $(-1)!!=1$.
			\hfill $\square$
		\end{corollary}
		
		\section{Appendix} 
		
		\begin{theorem}\label{thm:theorem6.1}
			For $k=1,2,3,...$, the following representation is valid:
			\[
			A_{2} \! ^{k} = e^{kz}\sum_{s=1}^{k}(-1)^{k-s}S(k,s)\left(\frac{d}{dz}\right)^{s}.
			\] 
		\end{theorem}
		\begin{proof}
			We shall prove the theorem using PMI.
			It is easy to see that there exist coefficients $b(k,s) \in \mathbb{Z}$ such that for $k=1, 2, 3, 4$
			\[
			A_{2}\!^{k} = e^{kz}\sum_{s=1}^{k}b(k,s)\left(\frac{d}{dz}\right)^{s}. 
			\]
			
			Let the above equality be true for some $k$. Then: \[
			A_{2}\!^{k+1} = A_{2} A_{2}\!^{k} = e^{(k+1)z} \left(kb(k,1)\frac{d}{dz} + \sum_{s=2}^{k}(kb(k,s) + b(k,s-1))\left(\frac{d}{dz}\right)^{s} + b(k,k)\left(\frac{d}{dz}\right)^{k+1}\right). 
			\]
			From the above equality we obtain the following three recurrence relations for $b(k,s)$:
			\begin{align*}
				b(k+1,1) &= kb(k,1) = k!;\\
				b(k+1,k+1) &= b(k,k) = 1;\\
				b(k+1,s) &= kb(k,s) + b(k,s-1),
			\end{align*}
			$s=2, \dots, k$.
			
			Let $b(k,s) = (-1)^{k+s}\widetilde{S}(k,s)$. Then, from the recurrence relations for b(k,s), we obtain the recurrence relations for the numbers $\widetilde{S}(k,s)$:
			\[
			\widetilde{S}(k,k) = 1; \; \widetilde{S}(k,1) = (-1)^{k-1}(k-1)!; \; \widetilde{S}(k+1,s) = \widetilde{S}(k,s-1) - k\widetilde{S}(k,s).
			\]
			
			But the Stirling numbers of the first kind $S(k,s)$ satisfy the same relations (see \cite{Abr}, p.~824). Thus $\widetilde{S}(k,s) = S(k,s)$ and the assertion of the theorem follows from PMI.
		\end{proof}
			
			\begin{theorem}\label{thm:theorem6.2}
				For $k=1,2,3,...$, the following representation is valid:
				
				\[
				A_{3} \! ^{k} = \sum_{s=0}^{k}(-1)^{k-s}(2k-2s-1)!!{2k-s-1 \choose s-1}\frac{1}{z^{2k-s}}\left(\frac{d}{dz}\right)^{3s},
				\]
				where we define $(-1)!! = 1$.
			\end{theorem}
			\begin{proof}
				Let $f \in U$. Then, 	with the help of PMI, one may establish that for $k\geq1$: 
				
				\[
				\bigg(\frac{1}{z}\frac{d}{dz}\bigg)^{k}f(z) = \sum_{s=1}^{k}\frac{a_{k}^{s}}{z^{2k-s}}f^{(s)}(z), 
				\]
				where $a_{k}^{s} \in \mathbb{Z}$ and  $a_{k}^{0} = 0$,  $a_{1}^{1}=1$, $a_{k}^{k}=1$. 
				
				Furthermore, we have 
				\begin{align*}
					\bigg(\frac{1}{z}\frac{d}{dz}\bigg)^{k+1}f(z) =& \  \bigg(\frac{1}{z}\frac{d}{dz}\bigg)\bigg[
					\bigg(\frac{1}{z}\frac{d}{dz}\bigg)^{k}f(z)\bigg] = \\[7pt]
					=& \ \frac{(-1)(2k-1)a_{k}^{1}}{z^{2(k+1)-1}}f^{(1)}(z) + \\[5pt] 
					& \ + \sum_{s=2}^{k}\frac{a_{k}^{s-1} - (2k-s)a_{k}^{s}}{z^{2(k+1) - s}}f^{(s)}(z) +\\[5pt]
					& \ + \frac{a_{k}^{k}}{z^{k+1}}f^{(s+1)}= \\[5pt]
					=& \sum_{s=1}^{k+1}\frac{a_{k+1}^{s}}{z^{2(k+1) - s}}f^{(s)}(z).
				\end{align*}
				
				Then the following relations hold: 		
				\begin{align}
					&1) \; a_{k+1}^{1} = (-1)(2k-1)a_{k}^{1}; \  a_{1}^{1} =1, \quad k=1,2,\dots;\nonumber\\
					&2) \; a_{k+1}^{k+1} = a_{k}^{k} = 1, \quad k=1,2,...;\label{47}\\
					&3) \; a_{k+1}^{s} = a_{k}^{s-1} - (2k-s)a_{k}^{s}, \quad
					s=2,3,...,k.\nonumber
				\end{align}
				
				With the help of 3) in (\ref{47}), for $s=k$, we obtain 
				\[
				a_{k+1}^{k} = a_{2}^{1} - (2+3+...+k) = -\frac{k(k+1)}{2}. 
				\]
				
				For $s=k-1$ we obtain 
				
				$$
				a_{k+1}^{k-1} = a_{3}-(4a_{3}^{2} + 5a_{4}^{3} + ... + (k+1)a_{k}^{k-1}).
				$$		
				i.e.	
				\[
				a_{k+1}^{k-1} = \frac{1}{2}\sum_{s=2}^{k}(s-1)s(s+1) = \frac{(k-1)k(k+1)(k+2)}{2.4}.
				\]
				For $s=k-2$ we obtain	
				\begin{align*}
					&a_{k+1}^{k-2} = a_{4}^{1} - \frac{1}{2}\cdot\frac{1}{4}\sum_{s=4}^{k}(s-2)(s)(s+1)(s+2) = \\
					& = -\frac{1}{2.4}\sum_{s=3}^{k}(s-2)(s-1)(s)(s+1)(s+2) = \\
					& = - \frac{(k-2)(k-1)k(k+1)(k+2)(k+3)}{2.4.6} \cdot
				\end{align*}
				
				The previous calculations lead us to the conjecture:
				
				\[
				a_{k+1}^{(k+1)-s} = (-1)^{s}\frac{(k-s+1)...(k-1)k...(k+s)}{2.4.6...(2s)} \cdot
				\]
				The last expression we rewrite in the form:
				
				\[
				a_{k+1}^{(k+1)-s} = (-1)^{s}(2s-1)!!{k+s \choose 2s}.
				\]
				After substituting $s$ with $k+1-s$ in the above equality, we obtain
				\[
				a_{k+1}^{s} = (-1)^{k+1-s}(2k-2s+1)!!{2k+1-s \choose s-1}, \ s=1,...,k; \ k=1,2,3,...
				\]
				
				Now, to prove strictly that the above representation of $a_{k+1}^{s}$ holds, it is enough to verify with its help that 
				1) and 3) are satisfied.

				The proof of 1) is trivial, while 3) follows after the substitution $2k-s = n$ and checking the equality
				
				\[
				(n+1-s){n+1 \choose n+2-s} = (n+1-s){n \choose n+2-s} + n{n-1 \choose n-s}, 
				\]
				with the help of the relation:
				
				\[
				(n+1-s){n \choose n+1-s} = n{n-1 \choose n-s}.
				\]
			\end{proof}
			
			\newpage 
			
				\section*{Acknowledgements} 
			
			The work has been (financially) supported by Grant KP-06-N 62/5
			of Bulgarian National Science Fund.

			\newpage
			
			\makeatletter
			\renewcommand{\@biblabel}[1]{[#1]\hfill}
			\makeatother


\begin{thebibliography}{99}
				
				
				\bibitem{Petkova} I. Petkova, \textit{Series expansion of the Gamma function and its reciprocal}, Notes on Number Theory and Discrete Mathematics, 27(4), 104-115, doi: 10.7546/nntdm.2021.27.4.104-115, 2021, \url{https://nntdm.net/papers/nntdm-27/NNTDM-27-4-104-115.pdf}
				
				\bibitem{Abr} Abramowitz, M. \& Stegun, I. A. (1972). \textit{Handbook of Mathematical Functions, Tenth Printing}, Washington. 
				
				
				
				\bibitem{7} Luke, Y. L. (1980). \textit{Special Mathematical Functions and Their Approximations}, Mir Publishing House, Moscow (in Russian).
				
				\bibitem{JEL} JANKE - EMDE - L$\ddot{O}$SH, \textit{Tafeln H$\ddot{o}$herer Funktionen}, Moscow, 1977, 
				(in Russian)
				
				\bibitem{Polya} G. Polya, G. Szeg\"o: \textit{Problems and Theorems in Analysis I}, Springer-Verlag Berlin Heidelberg, 1998 
				
				\bibitem{Ross} Clay C. Ross, Differential Equations, An Introduction with Mathematica (Second Edition), Springer, 2004
				
				\bibitem{Watson} G. N. Watson, A Treatise on the Theory of Bessel Functions (Second edition), Cambridge, England: Cambridge University Press, 1966. 
				
				\bibitem{Belbahri} Kamel Belbahri, Scale invariant operators and combinatorial expansions, Advances in Applied Mathematics 45 (2010) 548–563
				
				\bibitem{Riordan} J.   Riordan,   Combinatorial Identities, John Wiley \text{\&}  Sons, New York, 1968
				
				\bibitem{Kerimov} Kerimov, M.K. Studies on the zeros of Bessel functions and methods for their computation. Comput. Math. and Math. Phys. 54, 1337–1388 (2014).
				
				\bibitem{Vi} V. Kiryakova. Generalized Fractional Calculus and Applications (Pitman Research Notes in Mathematics Vol. 301, Longman, 1994)
				
				\bibitem{wiBell} \url{https://encyclopediaofmath.org/wiki/Bell_polynomial}
				
				\bibitem{eulerian} L.   Comtet,   Advanced Combinatorics , D. Reidel, Dordrecht, 1974
				
			\end{thebibliography}
		\end{document}